\documentclass[a4paper]{amsart}
\usepackage{graphicx}
\usepackage{amssymb}
\usepackage[mathscr]{eucal}
\usepackage{xspace}
\usepackage[latin1]{inputenc}
\usepackage{color}
\usepackage{hyperref}

\newcommand{\R}{\mathbb{R}}
\newcommand{\Rf}{{\mathfrak R}}

\newcommand{\Z}{\mathbb{Z}}

\newcommand{\LL}{\mathcal{L}}

\newcommand{\Oh}{\mathcal{O}}

\newcommand{\CS}{\mathcal{S}}

\newcommand{\te}{\theta}

\newcommand{\la}{\langle}
\newcommand{\ra}{\rangle}

\newcommand{\loc}{{\text{loc}}}
\newcommand{\X}{\times}

\renewcommand{\l}{\lambda}

\renewcommand{\a}{\alpha}
\renewcommand{\b}{\beta}

\newcommand{\s}{\sigma}
\newcommand{\g}{\gamma}
\renewcommand{\k}{\kappa}

\newcommand{\Dx}{\Delta x}
\newcommand{\Dt}{\Delta t}

\newcommand{\norm}[1]{\left\Vert#1\right\Vert}
\newcommand{\abs}[1]{\left\vert#1\right\vert}

\newcommand{\CC}{{\mathcal C}}

\DeclareMathOperator{\supp}{supp}

\newbox\bokstav
\newdimen\hoyde
\newcommand{\bgl}{{\hbox{$\left\lbrack\vbox to 8.5pt{}\right.\nOspace$}}}
\newcommand{\Bgl}{{\hbox{$\left\lbrack\vbox to 11.5pt{}\right.\nOspace$}}}
\newcommand{\bggl}{{\hbox{$\left\lbrack\vbox to 14.5pt{}\right.\nOspace$}}}
\newcommand{\Bggl}{{\hbox{$\left\lbrack\vbox to 17.5pt{}\right.\nOspace$}}}
\newcommand{\bgr}{{\hbox{$\left\rbrack\vbox to 8.5pt{}\right.\nOspace$}}}
\newcommand{\Bgr}{{\hbox{$\left\rbrack\vbox to 11.5pt{}\right.\nOspace$}}}
\newcommand{\bggr}{{\hbox{$\left\rbrack\vbox to 14.5pt{}\right.\nOspace$}}}
\newcommand{\Bggr}{{\hbox{$\left\rbrack\vbox to 17.5pt{}\right.\nOspace$}}}
\newcommand{\nOspace}{\nulldelimiterspace=0pt \mOth}
\newcommand{\mOth}{\mathsurround=0pt}
\newcommand{\Ljmp}{\mathopen{\lbrack\!\lbrack}}
\newcommand{\Rjmp}{\mathclose{\rbrack\!\rbrack}}
\newcommand{\bgLjmp}{\mathopen{\bgl\mskip-6mu\bgl}}
\newcommand{\bgRjmp}{\mathclose{\bgr\mskip-6mu\bgr}}
\newcommand{\BgLjmp}{\mathopen{\Bgl\!\!\Bgl}}
\newcommand{\BgRjmp}{\mathclose{\Bgr\!\!\Bgr}}
\newcommand{\bggLjmp}{\mathopen{\bggl\!\!\bggl}}
\newcommand{\bggRjmp}{\mathclose{\bggr\!\!\bggr}}
\newcommand{\BggLjmp}{\mathopen{\Bggl\!\!\Bggl}}
\newcommand{\BggRjmp}{\mathclose{\Bggr\!\!\Bggr}}
\newcommand{\jmp}[1]{%
\setbox\bokstav=\hbox{$ \left. #1\right. $}
\hoyde=\ht\bokstav 
\advance\hoyde by \dp\bokstav%
\hbox{$
        \ifinner
                \ifdim\hoyde<10pt
                   \Ljmp #1 \Rjmp%
                \else
                   \ifdim\hoyde <11pt
                      \Ljmp #1 \Rjmp%
                   \else
                      \ifdim\hoyde <14pt
                          \bgLjmp #1 \bgRjmp%
                      \else
                          \ifdim\hoyde <20pt
                             \BgLjmp #1 \BgRjmp%
                          \else
                              \bggLjmp #1 \bggRjmp%
                          \fi
                      \fi
                   \fi
                \fi
        \else
                \ifdim\hoyde<8.5pt
                   \Ljmp #1 \Rjmp%
                \else
                   \ifdim\hoyde <11.5pt
                      \bgLjmp #1 \bgRjmp%
                   \else
                      \ifdim\hoyde <14.5pt \Ch
                          \BgLjmp #1 \BgRjmp%
                      \else
                          \ifdim\hoyde <17.5pt
                             \bggLjmp #1 \bggRjmp%
                          \else
                              \BggLjmp #1 \BggRjmp%
                          \fi
                      \fi
                   \fi
            \fi
        \fi
$}
}

\theoremstyle{plain}
\newtheorem{theorem}{Theorem}[section]

\theoremstyle{definition}
\newtheorem{definition}{Definition}[section]
\theoremstyle{remark}

\newtheorem{remark}{Remark}[section]
\numberwithin{equation}{section}
\allowdisplaybreaks

\begin{document}

\title[$L^\infty$ solutions for polytropic gas with diffusive entropy]{$L^\infty$ solutions 
for a model of polytropic gas flow\\ with diffusive entropy}

\author[Frid]{Hermano Frid}
\address[Frid]{\newline
Institute of Pure and Applied Mathematics (IMPA), 
Est. Dona Castorina, 110, 
22460-320 Rio de Janeiro, Brazil} 

\author[Holden]{Helge Holden}
\address[Holden]{\newline
    Department of Mathematical Sciences,
    Norwegian University of Science and Technology,
    NO--7491 Trondheim, Norway,\newline
{\rm and} \newline
  Centre of Mathematics for Applications, 
University of Oslo,
  P.O.\ Box 1053, Blindern,
  NO--0316 Oslo, Norway }
\email[]{\href{mailto:holden@math.ntnu.no}{holden@math.ntnu.no}}
\urladdr{\href{http://www.math.ntnu.no/~holden}{www.math.ntnu.no/\~{}holden}}

\author[Karlsen]{Kenneth H.~Karlsen}
\address[Karlsen]{\newline
   Centre of Mathematics for Applications, 
 University of Oslo,
  P.O.\ Box 1053, Blindern,
  NO--0316 Oslo, Norway}
\email[]{\href{mailto:kennethk@math.uio.no}{kennethk@math.uio.no}}
\urladdr{\href{http://folk.uio.no/kennethk}{http://folk.uio.no/kennethk}}

\date{\today}

\thanks{Supported in part by the Research Council of Norway. This
 paper was written as part of  the international research program on
 Nonlinear Partial Differential Equations at the Centre for Advanced
 Study at the Norwegian Academy of Science and Letters in Oslo during
 the academic year 2008--09. H.~Frid gratefully acknowledges the 
 partial support from CNPq, grant 306137/2006-2, and 
 FAPERJ, grant E-26/152.192-2002.}

\begin{abstract} 
We establish the global existence of $L^\infty$ solutions for 
a model of polytropic gas flow with diffusive entropy. The result 
is obtained by showing the convergence of a class of finite difference 
schemes, which includes the Lax--Friedrichs and Godunov schemes.  
Such convergence is achieved by proving 
the estimates required for the application of the 
compensated compactness theory. 
\end{abstract}

\maketitle

\section{Introduction}
We consider the following  system modeling isentropic gas flow with smoothly varying entropy. 
The model reads in Eulerian coordinates
\begin{align}
& \rho_t +m_x=0, \label{e1}\\
& m_t+(\frac{m^2}{\rho}+p(\rho,S))_x=0,\label{e2}\\
& (\rho S)_t+(m S)_x=(\frac1{\rho}S_x)_x,\label{e3}   
\end{align}
where
$$
p(\rho,S)=\k e^{(\gamma-1)S/\Rf}\rho^{\gamma},
$$
where $\Rf>0$ and $\gamma>1$ are constants, and $\k=\frac{1}{4\g}(\g-1)^2$. 
Here $\rho$ represents the gas density, $m$ is the momentum defined as $m=\rho u$, where $u$ is the gas velocity, $p$ represents the gas pressure, and $S$ stands for the entropy. 
The system \eqref{e1}--\eqref{e3} is a mathematical model intended to approximate the more physical model where equation \eqref{e3} is replaced by the energy 
conservation law, which for smooth solutions is equivalent to the 
equation $(\rho S)_t+(m S)_x=0$, and this motivates our mathematical model. 

Initial data are given by
\begin{align}
& \rho(x,0)=\rho_0(x),& m(x,0)&=m_0(x), \label{e4}\\
&S(x,0)=S_0(x)=\s(y_0(x)),& y_0(x)&=\int_0^x\rho_0(z)\,dz.\label{e5}
\end{align}
Assume that
\begin{equation}
\rho_0,m_0,\frac{m_0}{\rho_0}\in L^\infty(\R),\qquad\rho_0\ge0,\qquad 
\s\in W_{\loc}^{3,2}(\R).\label{e6}
\end{equation}  
In particular, the initial data (and the solution) allows for the occurrence of vacuum. 
In addition, we also assume that $\s$ is periodic with period, say, $2\pi$, that is,
\begin{equation}\label{eper}
\s(y+2\pi)=\s(y),\qquad y\in\R.
\end{equation}

We remark that assumption \eqref{e6}, imposed on $\s$, implies that 
the solution of the heat equation with initial data $\s$,
\begin{equation}\label{e6'}
\tilde\s(y,t):=\frac{1}{(4\pi t)^{1/2}}\int_{\R} e^{-(y-z)^2/4t}\s(z)\,dz,
\end{equation}
satisfies
\begin{equation}\label{e6''}
|\tilde\s(y,t)-\bar\s|,|\tilde\s_y(y,t)|,|\tilde\s_{yy}(y,t)|\le C_0e^{-t},
\quad\text{with}\quad\bar\s:=\frac{1}{2\pi}\int_0^{2\pi}\s(z)\,dz,
\end{equation}
for some absolute constant $C_0>0$.
Indeed, \eqref{e6} and \eqref{eper} imply the absolute convergence 
of the Fourier series of $\s$, $\s'$ and $\s''$. On the other 
hand, a straightforward calculation shows that
\begin{align*}
\frac{1}{(4\pi t)^{1/2}}\int_{\R} e^{-(y-z)^2/4t} e^{ikz}\,dz&
=\frac{e^{(-y^2+(y+2ikt)^2)/4t}}{(4\pi t)^{1/2}}\int_{\R}e^{\frac{-(z-(y+2ikt))^2}{4t}}\,dz\\
&=e^{iky-k^2t},
\end{align*}
for any $k\in\R$, which then gives the asserted 
asymptotic behavior, by plugging the Fourier series 
for $\s$, $\s'$ and $\s''$ in \eqref{e6'} and the corresponding 
equations for $\tilde \s_y$ and $\tilde\s_{yy}$, obtained 
from \eqref{e6'} by replacing $\s$ by $\s'$ and $\s''$, respectively. 

We have the following definition of weak solution.
\begin{definition}\label{D:1} We say that  a function 
$(\rho,m,S)\in L^\infty(\R\X(0,\infty))$ is a weak solution to \eqref{e1}--\eqref{e5} if:
\begin{enumerate}
\item[(i)]  $m/\rho\in L^\infty(\R\X(0,\infty))$;
\item[(ii)] for all $\phi\in C_0^\infty(\R^2)$,
\begin{multline} 
\int_{\R\X(0,\infty)} (\rho,m)(x,t)\phi_t+(m,\frac{m^2}{\rho}+p(S,\rho))(x,t)\phi_x\,dx\,dt\\
 +\int_{\R}(\rho_0,m_0)(x)\phi(x,0)\,dx=0;\label{e8}
\end{multline}
\item[(iii)]  
\begin{align}
&S(x,t)=\frac1{\sqrt{4\pi t}}\int_{\R}e^{-\frac{(y(x,t)-z)^2}{4t}}\s(z)\,dz,\label{e9}\\
\intertext{where}
&y(x,t)=\int_0^x\rho(z,t)\,dz-\int_0^tm(0,s)\,ds.\label{e10}
\end{align}
\end{enumerate}
\end{definition}

We observe that away from vacuum, equation \eqref{e3}, through 
the Lagrange transformation, $(x,t)\mapsto(y(x,t),t)$, with 
$y(x,t)$ given by \eqref{e10}, becomes
$$
S_t=S_{yy},
$$
and this justifies (iii) of Definition~\ref{D:1}.

Indeed, the interplay between the Lagrangian and Eulerian 
formulation of the model is important. For the record 
we note that the model \eqref{e1}--\eqref{e3} reads 
in Lagrangian coordinates
\begin{equation}\label{eq:lagrange}
\begin{aligned}
v_t-u_y&=0, \\
u_t+p(v,S)_y&=0, \\
S_t&=S_{yy}
\end{aligned}
\end{equation}
where $v=1/\rho$ is the specific volume. We remark that, despite 
the fact that system \eqref{eq:lagrange} has a 
form much simpler than \eqref{e1},\eqref{e2},\eqref{e3}, the 
possibility of occurrence of vacuum turns the direct analysis of the 
Cauchy problem for \eqref{eq:lagrange} a 
very difficult task and so, as in the isentropic case, a better 
strategy is to proceed with the the analysis of 
the corresponding problem in 
Eulerian coordinates, that is,  \eqref{e1}--\eqref{e5}.

Our main result reads as follows.
\begin{theorem}\label{T:1}  
There exists a constant $r(\g)>0$ such that if $\|(\rho_0,m_0)\|_\infty<r(\g)$, then 
there exists  a global weak solution to the Cauchy 
problem \eqref{e1}--\eqref{e5} satisfying an entropy inequality of the form
\begin{equation}\label{e100}
\eta_*(\rho,m,S)_t+q_*(\rho,m,S)_x\le -Ce^{-t},
\end{equation}
in the sense of distributions, for some $C>0$ depending on $L^\infty$ bounds for $\rho,m,S$, where  
\begin{equation}\label{eent}
\eta_*(\rho,m,S)= \frac12\rho u^2+\frac{\k}{\g-1}e^{(\g-1)S/\Rf}\rho^\g,\qquad 
q_*(\rho,m,S)=u\eta_*(\rho,m,S)+pu.
\end{equation}
Moreover, $r(\g)\to\infty$ as $\g\to1+$. Further, if $\rho_0,m_0$ are 
periodic with period $L$ such that  $y_0(L)=2\pi$, we have the following decay
\begin{equation}\label{edecay}
\lim_{t\to\infty}\int_0^L|\left(\rho(x,t),m(x,t),S(x,t)\right)-(\bar\rho,\bar m,\bar S)|\,dx=0,
\end{equation}    
where $\bar\rho,\bar m,\bar S$ are the mean values of $\rho_0,m_0,S_0$, respectively. 
\end{theorem}

\section{Background results}
Let us first recall results for the $p$-system for a polytropic 
gas in Eulerian coordinates. More precisely, we consider the system
\begin{align}
 \rho_t +m_x&=0, \label{e1A}\\
 m_t+\big(\frac{m^2}{\rho}+p(\rho)\big)_x&=0,\label{e2B}
\end{align}
where the pressure is given by $p(\rho)=\k e^{(\gamma-1)S/\Rf}\rho^{\gamma}$. 
For later use we observe that we can rewrite the conserved 
quantities in terms of the other variables, viz.,
\begin{equation} \label{eq:reskrive}
\rho=\rho(p,S)=\big(\frac{p}{\k}\big)^{1/\g}e^{(\g-1)S/(\g \Rf)}, \quad
m=m(u,p,S)=\rho u= u\big(\frac{p}{\k}\big)^{1/\g}e^{(\g-1)S/(\g \Rf)}.
\end{equation}
Here we consider the isentropic case where the entropy $S$ is considered a constant. 
Recall that the functions
\begin{align}
w&=u+\frac{1}{\te}\big(p_\rho\big)^{1/2}=u+e^{\te S/\Rf}\rho^{\te}
=u+\big(\frac{p}{\k}\big)^{\te/\g}e^{-\te S/(\g \Rf)},\\ 
z&=u-\frac{1}{\te}\big(p_\rho\big)^{1/2}=u-e^{\te S/\Rf}\rho^{\te}
=u-\big(\frac{p}{\k}\big)^{\te/\g}e^{-\te S/(\g \Rf)},
\end{align}
with $\te=\frac12(\g-1)$,
form a pair of Riemann invariants for system \eqref{e1A}--\eqref{e2B} in 
the isentropic case where $S$ is constant. A standard 
calculation (see, e.g., \cite{HR,DCL}) yields that the rarefaction curves are given by
$$
m=\frac{m_l}{\rho_l}\rho\pm \g^{1/2}e^{\te S/\Rf}\rho(\rho^\te-\rho_l^\te),
$$
while the Hugoniot locus reads
$$
m=\frac{m_l}{\rho_l}\rho\pm\te e^{\te S/\Rf}\rho  
\Big(\frac{1}{\rho\rho_l}(\rho^\g-\rho_l^\g)(\rho-\rho_l)\Big)^{1/2},
$$
from a given left state $(\rho_l,m_l)$.
When we involve the entropy condition we find that the wave curves equal
\begin{align}\label{eq:W1}
&W_1(\rho_l,m_l): \quad m= \frac{m_l}{\rho_l}\rho- 
\begin{cases} \g^{1/2}e^{\te S/\Rf}\rho(\rho^\te-\rho_l^\te) & \text{for $\rho\le \rho_l$},\\
\te e^{\te S/\Rf}\rho  \Big(\frac{1}{\rho\rho_l}(\rho^\g-\rho_l^\g)(\rho-\rho_l)\Big)^{1/2}
& \text{for $\rho\ge \rho_l$},
 \end{cases}\\[2mm]
&W_2(\rho_l,m_l):\quad m= \frac{m_l}{\rho_l}\rho+ \begin{cases}
\te e^{\te S/\Rf}\rho  \Big(\frac{1}{\rho\rho_l}(\rho^\g-\rho_l^\g)(\rho-\rho_l)\Big)^{1/2} 
& \text{for $\rho\le \rho_l$},\\
 \g^{1/2}e^{\te S/\Rf}\rho(\rho^\te-\rho_l^\te)& \text{for $\rho\ge \rho_l$}.
 \end{cases}\label{eq:W2}
\end{align}
In the variables $(\rho,u)$ we find
\begin{align}\label{eq:W1u}
&W_1(\rho_l,u_l): \quad u= u_l- \begin{cases} \g^{1/2}e^{\te S/\Rf}(\rho^\te-\rho_l^\te) & \text{for $\rho\le \rho_l$},\\
\te e^{\te S/\Rf}\Big(\frac{1}{\rho\rho_l}(\rho^\g-\rho_l^\g)(\rho-\rho_l)\Big)^{1/2}& \text{for $\rho\ge \rho_l$},
 \end{cases}\\[2mm]
&W_2(\rho_l,u_l):\quad u= u_l+ \begin{cases}\te e^{\te S/\Rf}  \Big(\frac{1}{\rho\rho_l}(\rho^\g-\rho_l^\g)(\rho-\rho_l)\Big)^{1/2} & \text{for $\rho\le \rho_l$},\\
 \g^{1/2}e^{\te S/\Rf}(\rho^\te-\rho_l^\te)& \text{for $\rho\ge \rho_l$}.
 \end{cases}\label{eq:W2u}
\end{align}
\begin{figure}[tbp]
  \centering
  \includegraphics[width=0.45\linewidth]{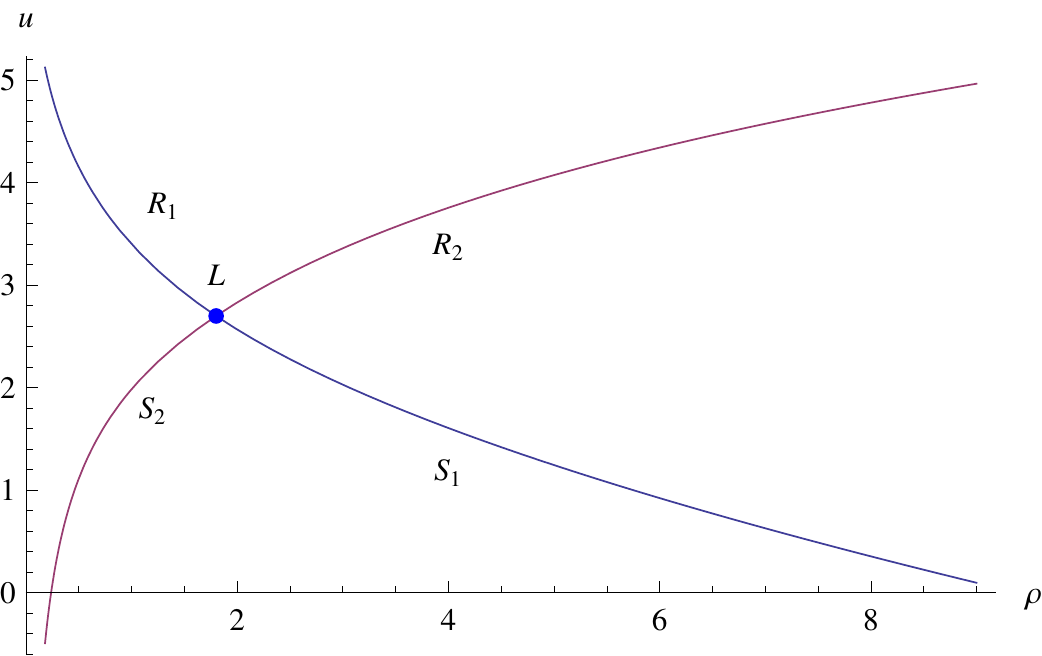} 
  \includegraphics[width=0.45\linewidth]{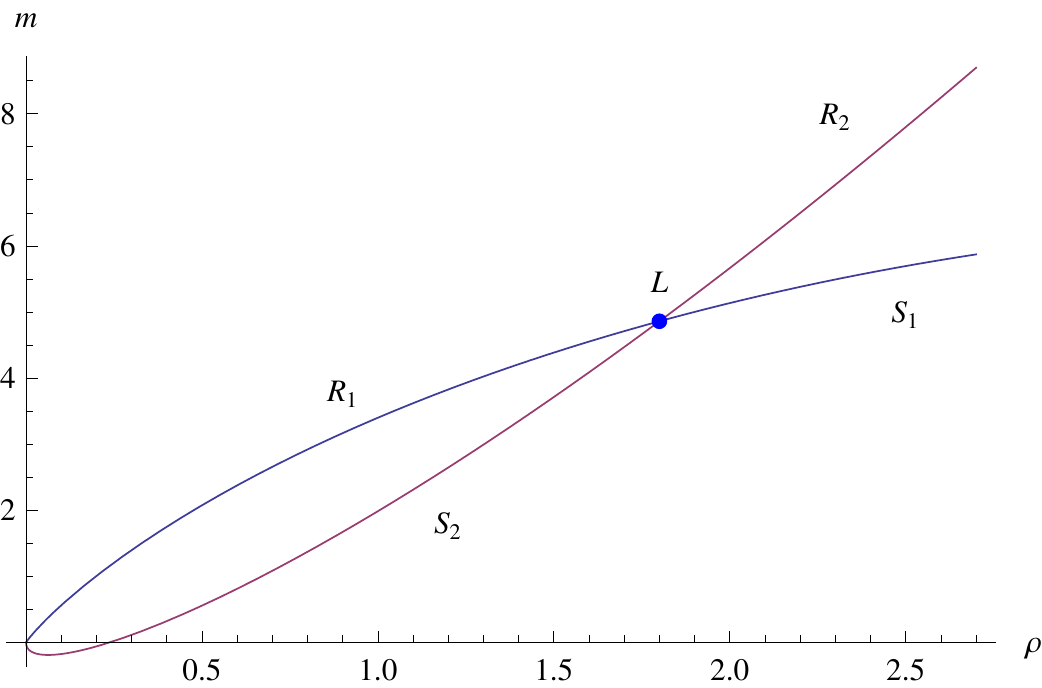}
  \caption{The wave curves for the $p$-system for a given left state.} 
  \label{fig:wave}
\end{figure}

An important property of the $p$-system is that the Riemann invariants provide  invariant regions. More specifically, (see, e.g., \cite[Lemma 5]{DCL}) if $(\rho_0(x),m_0(x))\in \Omega=\{(\rho,m)\mid w\le w_0, \, z\ge z_0, \, w-z\ge 0\}$ for all $x\in\R$, then also the solution $(\rho(x,t),m(x,t))$ will remain in $\Omega$, that is, $(\rho(x,t),m(x,t))\in \Omega$ for $(x,t)\in\R\X[0,\infty)$.

\begin{figure}[tbp]
  \centering
  \includegraphics[width=0.5\linewidth]{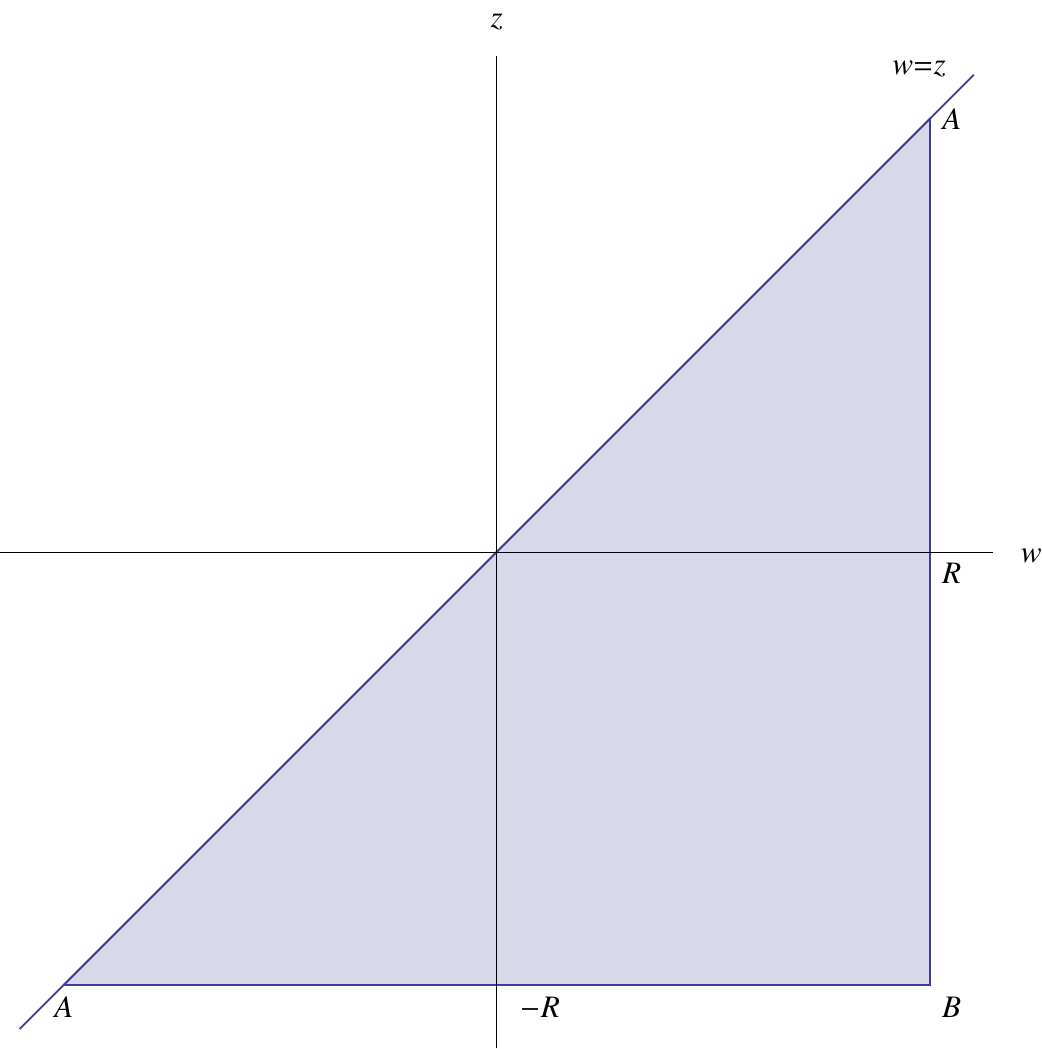}  \includegraphics[width=0.45\linewidth]{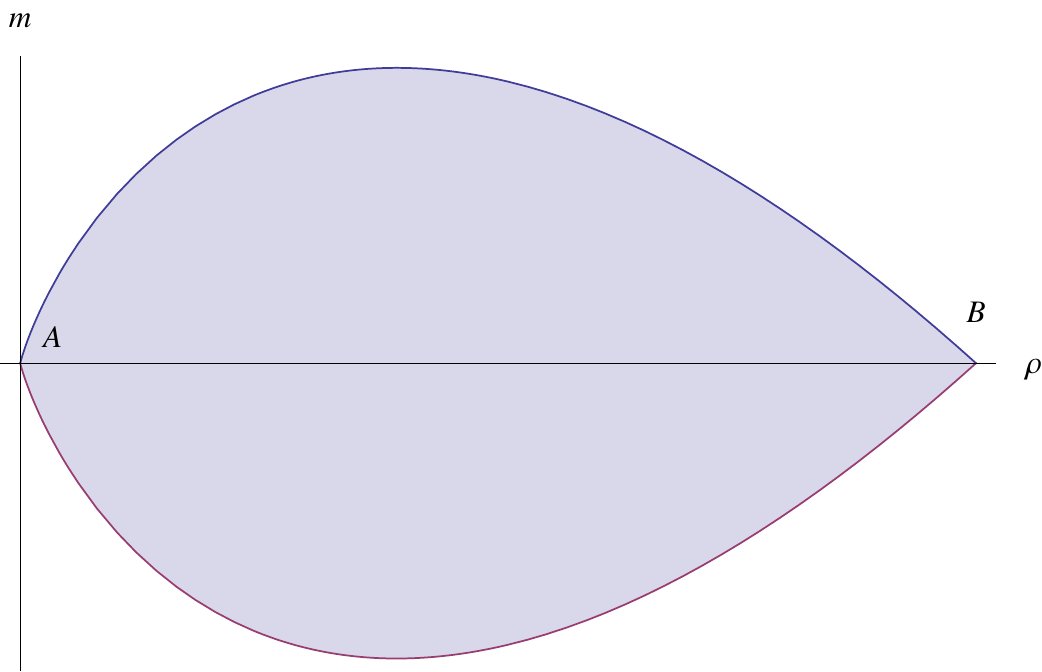}
  \caption{The triangle in the Riemann invariants (left) is mapped into the indicated region bounded by rarefaction waves in the 
  $(\rho,m)$-plane.} 
  \label{fig:trekant}
\end{figure}
\begin{figure}[tbp]
  \centering
  \includegraphics[width=0.4\linewidth]{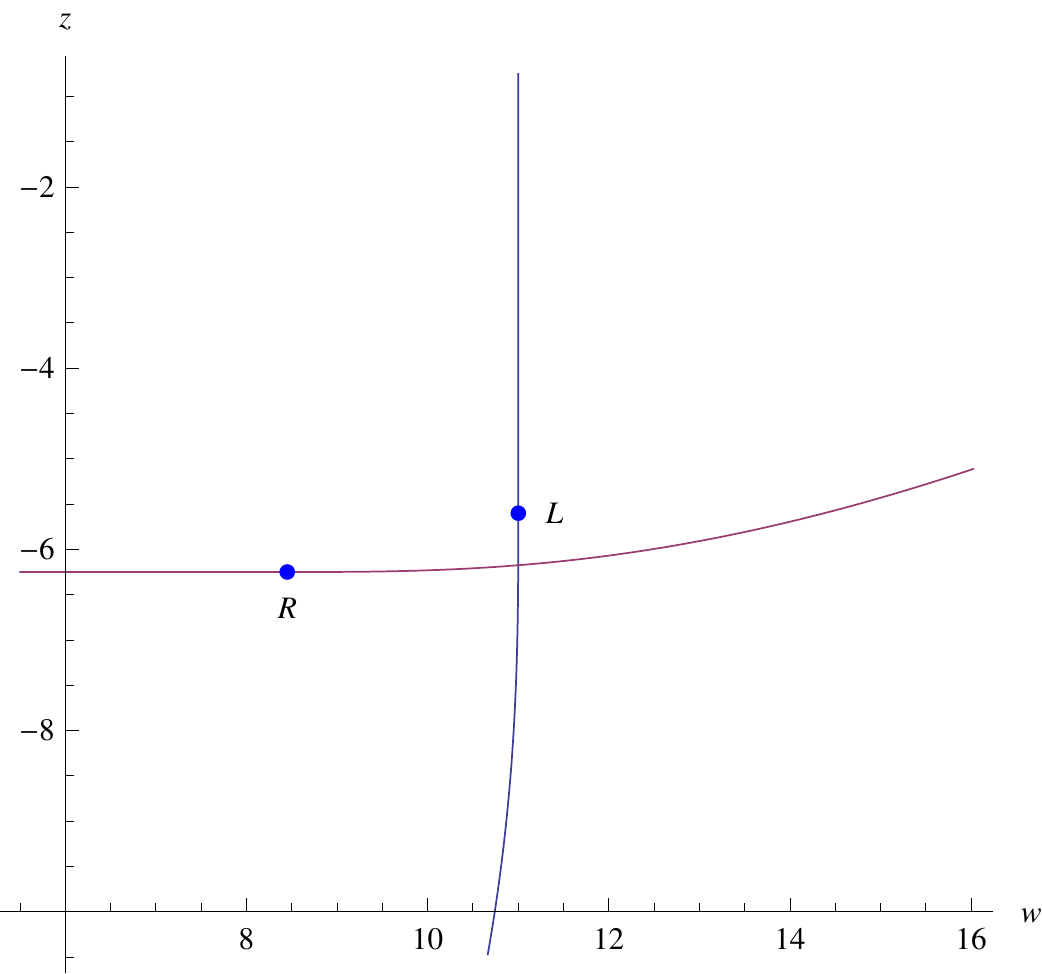}
  \caption{The solution of the Riemann problem using 
  the Riemann invariants as coordinates. Through the right state 
  backward Riemann invariants are drawn.} 
  \label{fig:invR}
\end{figure}

An entropy-entropy flux pair $(\eta, q)$ for the $p$-system satisfies for smooth solutions 
$$
\eta(\rho,m)_t+q(\rho,m)_x=0.
$$
Consistency with the system \eqref{e1A}--\eqref{e2B} requires
\begin{equation}\label{eq:consist}
\nabla q(\rho,m)=\nabla \eta(\rho,m)\nabla F(\rho,m),
\end{equation}
where $F=(m, \frac{m^2}{\rho}+p(\rho))$ is the flux function of the $p$-system. 
A particular choice of entropy-entropy flux pair $(\eta_*, q_*)$ reads
\begin{align}
\eta_*&=\frac{m^2}{2\rho}+\frac{p}{\g-1}= \frac12\rho u^2
+\frac{\k}{\g-1}e^{(\g-1)S/R}\rho^\g,\\ 
q_*&=u\eta_*+ pu=u\eta_*+ u\k e^{(\g-1)S/R}\rho^\g.
\end{align}

More generally, the weak entropy-entropy flux pairs  $(\eta,q)$ constitute 
a class of entropy-entropy flux pairs of particular interest in isentropic gas 
dynamics, as first pointed out in \cite{DP2}, and they are 
characterized by the following conditions at the vacuum line:
$$
\eta(\rho,u)|_{\rho=0}=0,\qquad \eta_\rho(\rho,u)|_{\rho=0}=g(u),
$$
for some continuous function $g$. Let us denote
$$
\chi(\rho,u;S)=(\frac{p}{\rho}-u^2)_+^\l,
$$
where $(a)_+=\max\{0,a\}$  and $\l=\frac{3-\g}{2(\g-1)}$. 
As observed in \cite{LPT}, weak entropy-entropy flux pairs 
can be given by the integral formulas
\begin{align}
& \eta(\rho,u)=\int_{\R}g(\xi)\,\chi(\rho,\xi-u)\,d\xi,\label{ewep1}\\
&q(\rho,u)=\int_{\R}g(\xi)\big(\theta\xi+(1-\theta)u\big)\chi(\rho,\xi-u)\,d\xi.\label{ewep2}
\end{align}

\begin{remark}\label{R:1} Observe that the entropy pair $(\eta_*,q_*)$, 
defined in \eqref{eent}, is a weak convex entropy pair. 
Moreover, for any weak entropy pair $(\eta,q)$ there 
exists a constant $C_\eta>0$ such that $\eta+C_\eta \eta_*$ is convex. 
\end{remark}

\begin{figure}[tbp]
  \centering
  \includegraphics[width=0.7\linewidth]{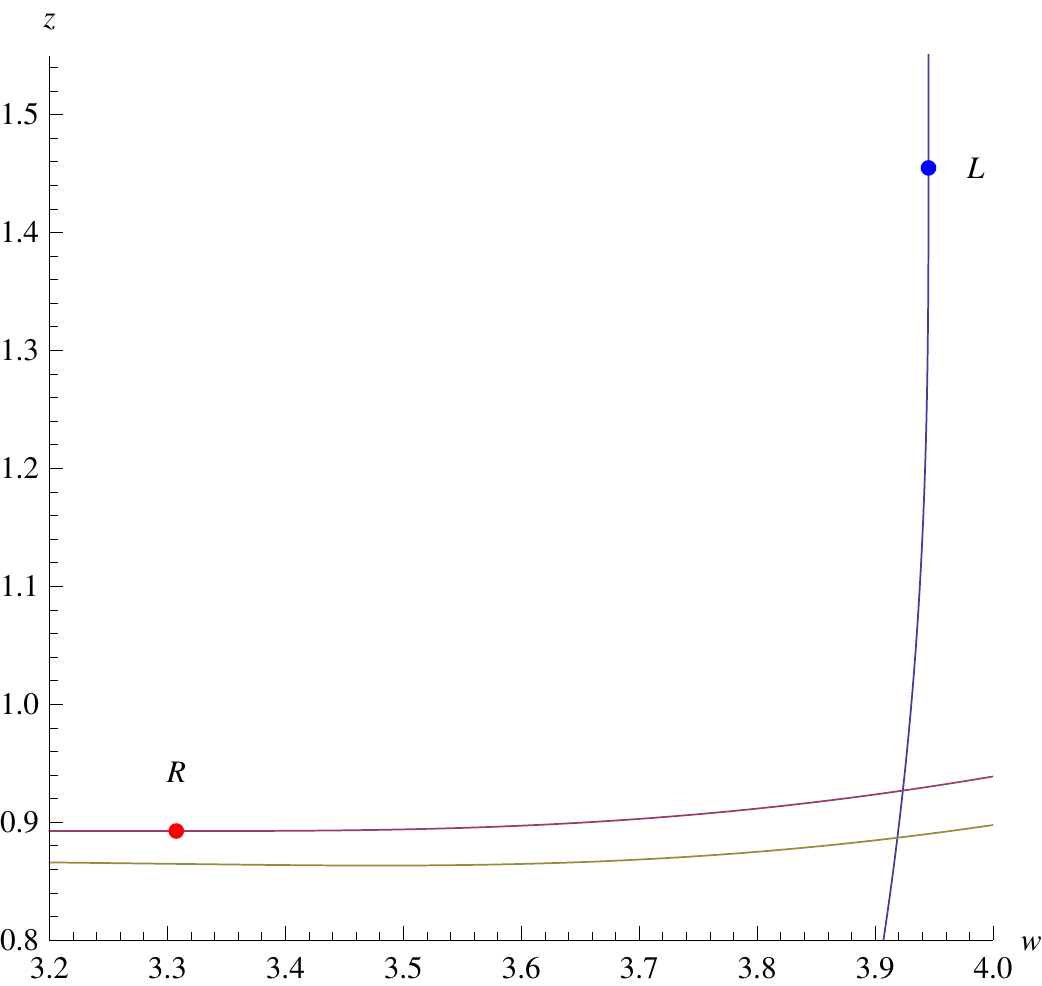}
  \caption{The slow Riemann invariant through the left state (blue curve), 
  and the backward fast Riemann invariant through the right state (red curve). 
  In addition the yellow curve $(\bar w,\bar z)$, whose intersection with the slow Riemann 
  invariant determines the contact discontinuity.} 
  \label{fig:RP}
\end{figure}

\begin{figure}[tbp]
  \centering
  \includegraphics[width=0.45\linewidth]{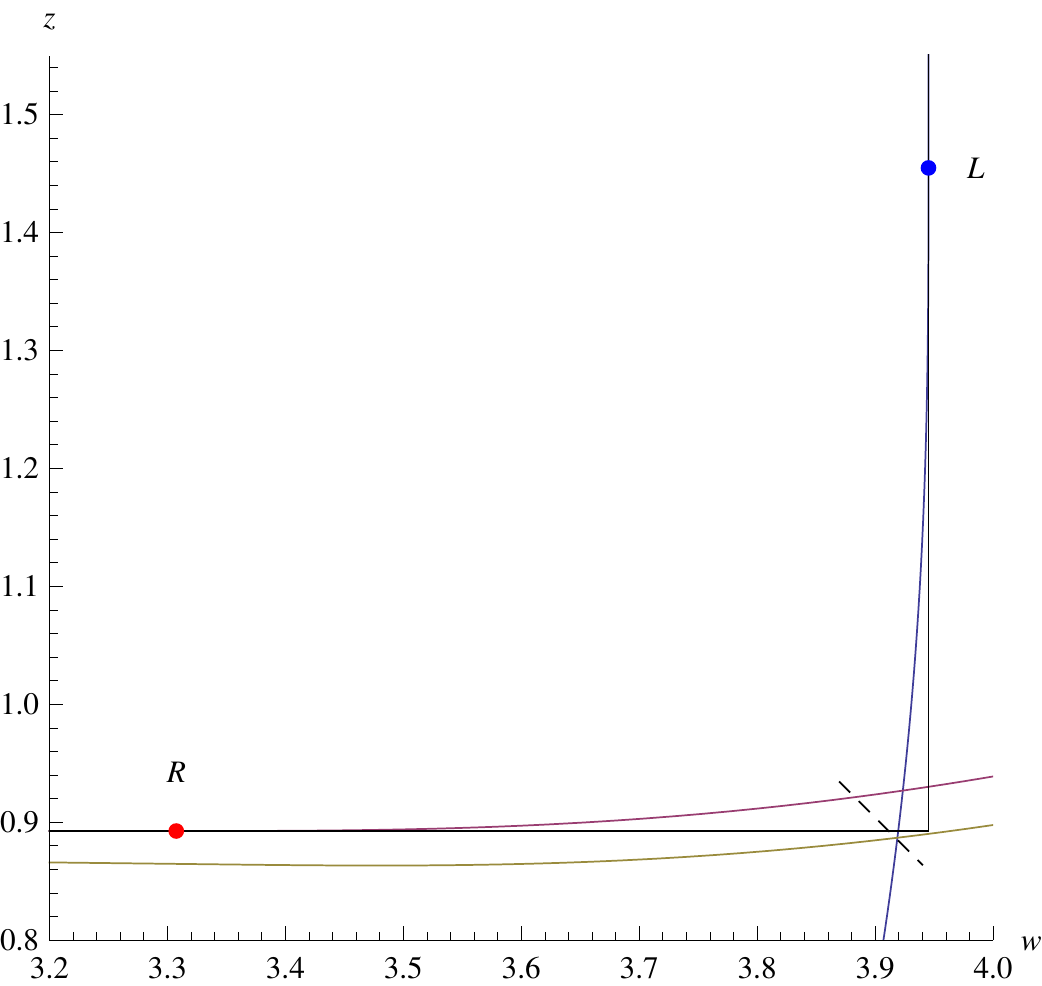}
   \includegraphics[width=0.45\linewidth]{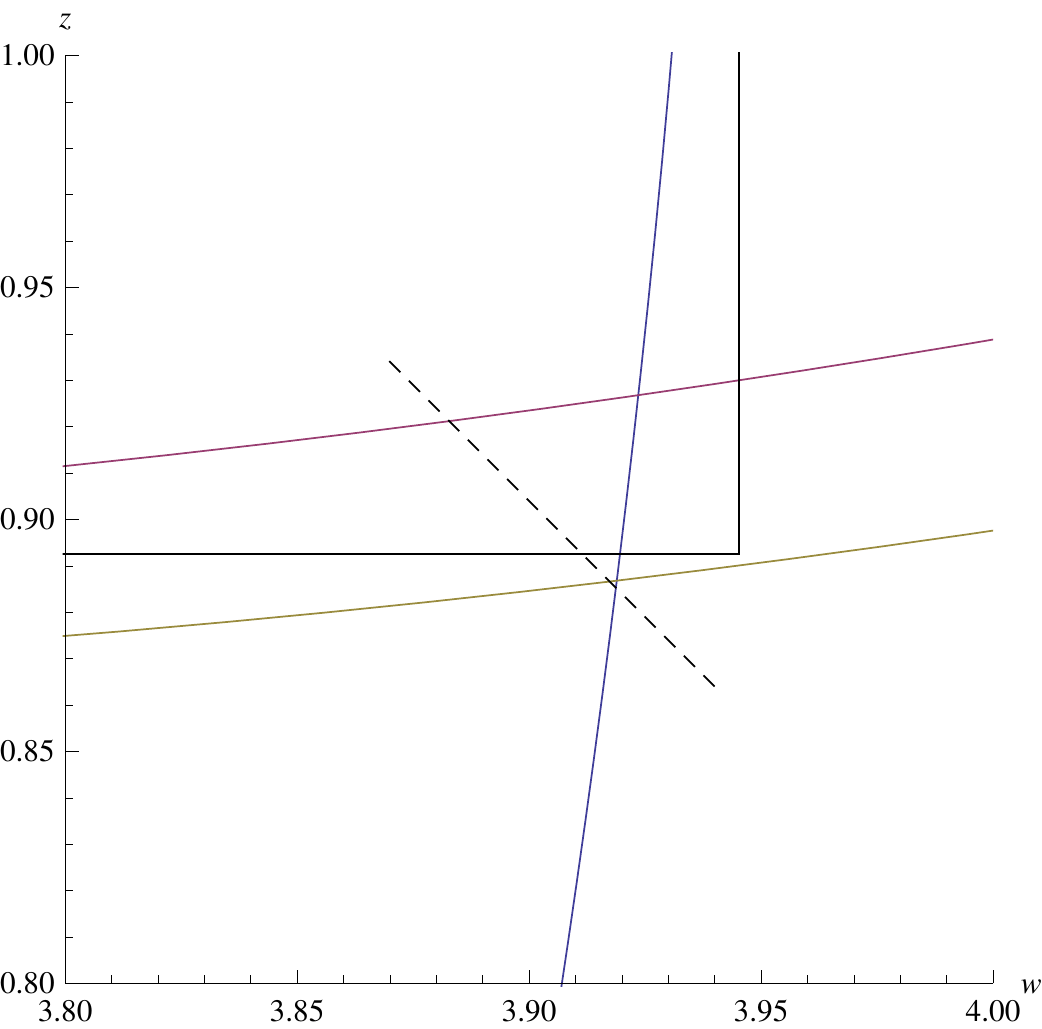}
  \caption{The same data is in Figure \ref{fig:RP}. Curves for the invariant region 
  for the corresponding $p$-system are added (black). In addition, the 
  dashed line is given by $w-z$ equals 
  a constant determined by the intersection between the yellow and blue curves. 
  The interaction of this straight line with the red curve gives the value 
  on the right of the contact discontinuity.  
  The right figure is a close-up near the intersection.} 
  \label{fig:RP2a}
\end{figure}

\medskip
Let us now turn to the full system
\begin{align}
& \rho_t +m_x=0, \label{e1'a}\\
& m_t+(\frac{m^2}{\rho}+p(\rho,S))_x=0,\label{e2'b}\\
& (\rho S)_t+(m S)_x=0,\label{e3'c}   
\end{align}
where the pressure $p$ is given as above.
The Riemann problem is the initial value problem for the 
system \eqref{e1'a}--\eqref{e3'c} with   special initial data consisting 
of a single jump between two constant states, viz.
$$
\begin{pmatrix} \rho \\ m \\ \rho S  \end{pmatrix}\Bigg|_{t=0}(x)=
\begin{cases}
\begin{pmatrix} \rho_l \\ m_l \\ \rho_l S_l  \end{pmatrix}& \text{for $x<0$},\\[8mm]
\begin{pmatrix} \rho_r \\ m_r \\ \rho_r S_r  \end{pmatrix} & \text{for $x>0$}.
\end{cases}
$$
The system \eqref{e1'a}--\eqref{e3'c} possesses three eigenfields associated with the eigenvalues 
$$
\l_1=u-\sqrt{p_\rho},\qquad \l_2=u,\qquad \l_3=u+\sqrt{p_\rho}.
$$
The solution to a Riemann problem for system \eqref{e1'a}--\eqref{e3'c} may be described using the coordinates $w,z,S$, that is, the Riemann invariants for the $p$-system and the entropy, in the following way. Consider first the case when the solution does not contain vacuum.  
The solution of the Riemann problem, starting from the left state $(\rho_l,m_l,S_l)$, consists of a slow wave in which the entropy $S$ remains constant (i.e., in the $(w,z)$-plane determined by $S=S_l$), followed by a contact discontinuity in which the velocity $u$ and the pressure $p$ remain unchanged, and finally a fast wave with constant entropy $S$ (i.e., in the $(w,z)$-plane determined by $S=S_r$) connected with the given right state 
$(\rho_r,m_r,S_r)$.  Along the slow wave we can write the Riemann invariants as\footnote{It turns out to be easier to describe the solution using the speed $u$ rather than the momentum $m$ as a variable.}
\begin{equation}\label{eq:wz1}
\begin{aligned}
w&=u_1(\rho;\rho_l,u_l,S_l)+e^{\te S_l/\Rf}\rho^\te, \\
z&=u_1(\rho;\rho_l,u_l,S_l)-e^{\te S_l/\Rf}\rho^\te,
\end{aligned}
\end{equation}
where $u=u_1(\rho;\rho_l,u_l,S_l)$ is the slow wave given by \eqref{eq:W1}.   
For the fast wave we consider the \textit{backward} wave (i.e., consisting of the states that can be connected to a given right state from the left), and the Riemann invariants read
\begin{equation}\label{eq:wz2}
\begin{aligned}
w&=\tilde u_2(\rho;\rho_r,u_r,S_r)+e^{\te S_r/\Rf}\rho^\te, \\
z&=\tilde u_2(\rho;\rho_r,u_r,S_r)-e^{\te S_r/\Rf}\rho^\te,
\end{aligned}
\end{equation}
where $u=\tilde u_2(\rho;\rho_r,u_r,S_r)$ is the fast  
\textit{backward} wave corresponding to \eqref{eq:W2}. 
The contact discontinuity, with pressure $p^*$ and velocity $u^*$, jumps 
from a left density $\rho_l^*$ to a right density   
$\rho_r^*$ determined by 
\begin{equation}\label{eq:wz1/2}
\begin{aligned}
p^*&=\kappa e^{(\g-1)S_l/\Rf} (\rho_l^*)^\g=\kappa e^{(\g-1)S_r/\Rf} (\rho_r^*)^\g, \\
u^*&=u_1(\rho_l^*;\rho_l,u_l,S_l)=\tilde u_2(\rho_r^*;\rho_r,u_r,S_r),
\end{aligned}
\end{equation}
which yields
\begin{equation}\label{eq:wz1/2a}
\frac{\rho_l^*}{\rho_r^*}= e^{(\g-1)(S_r-S_l)/(\g \Rf)}
\end{equation}
to be inserted in the second equation for the velocity, $u_1=\tilde u_2$, to determine 
$\rho_l^*$ and $\rho_r^*$. In terms of the Riemann invariants we find that $w$ jumps from 
$u^*+(p^*/\k)^{\te/\g} e^{-\te S_l/(\g \Rf)}$ to $u^*+(p^*/\k)^{\te/\g} e^{-\te S_r/(\g \Rf)}$, and 
similarly  $z$ jumps from 
$u^*-(p^*/\k)^{\te/\g} e^{-\te S_l/(\g \Rf)}$ to $u^*-(p^*/\k)^{\te/\g} e^{-\te S_r/(\g \Rf)}$.  
An alternative way to describe the contact discontinuity is the following. 
Consider a point on the backward fast wave curve with 
Riemann invariants $(w,z)$ given by \eqref{eq:wz2}, which we can write as
$w=\tilde u_2+(p/\k)^{\te/\g}e^{-\te S_r/(\g \Rf)}$ and 
$z=\tilde u_2-(p/\k)^{\te/\g}e^{-\te S_r/(\g \Rf)}$. 
Construct now another curve $(\bar w,\bar z)$, given 
as a Riemann invariant with the same velocity $\tilde u_2$ and 
pressure $p$ as $(w,z)$, but with the entropy $S_r$ replaced by $S_l$, that is,
$$
\bar w=\tilde u_2+(p/\k)^{\te/\g}e^{-\te S_l/(\g \Rf)}, \quad 
\bar z=\tilde u_2-(p/\k)^{\te/\g}e^{-\te S_l/(\g \Rf)}.
$$
We find
\begin{align*}
w+z&=2\tilde u_2=\bar w+\bar z, \\
w-z&=2(p/\k)^{\te/\g}e^{-\te S_r/(\g \Rf)}=(\bar w-\bar z) e^{\te (S_l-S_r)/(\g \Rf)},
\end{align*}
which yields
\begin{align*}
\bar w&=\frac{w}{2}(1+e^{\te (S_l-S_r)/(\g \Rf)})+\frac{z}{2}(1-e^{\te (S_l-S_r)/(\g \Rf)}),\\ 
\bar z&=\frac{w}{2}(1-e^{\te (S_l-S_r)/(\g \Rf)})+\frac{z}{2}(1+e^{\te (S_l-S_r)/(\g \Rf)}).
\end{align*}
The intersection between the slow wave curve in the  Riemann 
invariants plane and the curve $(\bar w,\bar z)$ determines the values 
of the variables to the left of the contact discontinuity. Through this 
intersection we draw the line where $w+z$ is constant, and 
the intersection between this line and the backward fast 
wave gives the values of the variables to the 
right of the contact discontinuity,  cf.~Figures \ref{fig:RP} and  \ref{fig:RP2a}.


The solution involves vacuum when the slow wave is a rarefaction 
wave that connects to a state on the vacuum line $w=z$; the velocity 
is then given by $u^*=u_l+\g^{1/2}e^{\te S_l/\Rf}\rho_l^\te$ and $w=z=u^*$.  
Similarly, the given right state connects via a rarefaction 
from a vacuum state with velocity $\tilde u^*=u_r-\g^{1/2}
e^{\te S_r/\Rf}\rho_r^\te$ and $w=z=\tilde u^*$.

\section{Proof of Theorem \ref{T:1}}

\subsection{Construction of approximate solutions}
Here we provide the full proof of Theorem \ref{T:1}.  
We construct approximate solutions for \eqref{e1}--\eqref{e3} by using a 
Godunov-type finite difference scheme based on 
solving Riemann problems at each time step, updating the 
approximate $S$ using the Lagrange 
transformation, and averaging at the end of each time step. 

Before we begin the proof, let us describe the fundamentals of 
the construction of the approximate solution. We discretize both 
in space and time.  Let $h=\Delta t$, and 
$\Delta x=ch$ with $c>0$ to be chosen by the CFL condition 
$$
c>\sup_{(x,t)\in\R\X[0,\infty)}\left|\frac{m^h(x,t)}{\rho^h(x,t)}\pm 
\sqrt{p_\rho(\rho^h(x,t),S^h(x,t))}\right|,
$$
which is possible as long as we can obtain an $L^\infty$ a priori 
bound for 
$$
\frac{m^h(x,t)}{\rho^h(x,t)}\pm \sqrt{p_\rho(\rho^h(x,t),S^h(x,t))}.
$$ 

The initial data $\rho_0, m_0, S_0$ is approximated by step functions with 
jumps at $x_{i-1/2}:=(i-1/2)\Delta x$ for $i\in\Z$.  The multiple Riemann 
problems are solved for $t\in[0,h)$. At $t=h$ a new step function is created 
with jumps at $x_{i-1/2}$ (details given below), and new Riemann 
problems are solved. More precisely, suppose the 
approximate solution $U^h=(\rho^h,m^h,\rho^h S^h)$ has 
been defined for $t\le jh$ and that 
$U^h(x,jh)$ is constant for $x\in I_i$ where
$$
I_i=(x_{i-1/2},x_{i+1/2}), \quad i\in\Z.
$$
For $t\in[jh,(j+1)h)$, setting $x_i=i\Delta x$, $i\in\Z$,  we define $U^h(x,t)$ by 
glueing together the solutions of the Riemann 
problems for the system  \eqref{e1'a}--\eqref{e3'c}
defined at $[x_i,x_{i+1}]\X[jh,(j+1)h)$, determined by the discontinuities 
at the points $(x_{i+1/2},jh)$, $i\in\Z$.  Inductively this yields a function 
$U^h$ defined on $\R\X[0,\infty)$, as long as we are able to obtain 
the necessary  a priori bound mentioned above.   


We now provide the details of the construction of the 
approximate solution. Assume that we have constructed the 
approximate solution $U^h$ for $x\in\R$ and $t<jh$, and 
have defined it at time $t=jh$ as a piecewise constant function 
with jumps at $x_{i+1/2}$ for $i\in\Z$. For $(x,t)\in [x_i,x_{i+1}]\X[jh,(j+1)h)$, 
$i\in\Z$, let $U^h(x,t)$  be the solution of the 
Riemann problem \eqref{e1'a}--\eqref{e3'c} as described in the previous section.   Set
\begin{align*}
y^h(x,t)&=\int_0^x\rho^h(z,t)\,dz-\int_0^t m^h(0,s)\,ds,
\quad x\in\R,\ t\in[jh,(j+1)h),\\
\intertext{and}
\s^h(x,t)&=\frac1{\sqrt{4\pi t}}\int_{\R}e^{-\frac{(y^h(x,t)-z)^2}{4t}}\s(z)\,dz
=\tilde\s(y^h(x,t),t),\quad (x,t)\in\R\X[0,(j+1)h).
\end{align*}
We then define\footnote{We use the standard notation 
$f(x\pm 0)=\lim_{\epsilon \downarrow 0}f(x\pm\epsilon)$.} 
\begin{align}
 \rho^h(x,(j+1)h)&=\frac{1}{\Delta x}\int_{I_i} \rho^h(\tilde x,(j+1)h-0)\,d \tilde x,\label{eq:snittB}\\
 m^h(x,(j+1)h)&=\frac{1}{\Delta x}\int_{I_i}  m^h(\tilde x,(j+1)h-0)\,d \tilde x,\label{eq:snittB2}\\
S^h(x,(j+1)h)&=\frac{1}{\Delta x}\int_{I_i}\s^h(\tilde x,(j+1)h-0)\,d \tilde x,\label{eq:snitt}
\end{align}
for $x\in I_i$.

\subsection{Convergence proof}
We now address the questions of the $L^\infty$ a priori bound 
and convergence of $U^h$ as $h\to 0$. First, we investigate the 
problem of obtaining an a priori $L^\infty$ bound 
for the approximate solution. Let us denote
$$
w^h(x,t)=w(U^h(x,t)),\qquad z^h(x,t)=z(U^h(x,t)).
$$
Let $r>0$ be such that
$$
w^h(x,0)\le r,\qquad z^h(x,0)\ge -r,\qquad x\in\R.
$$
\begin{figure}[tbp]
  \centering
  \includegraphics[width=0.5\linewidth]{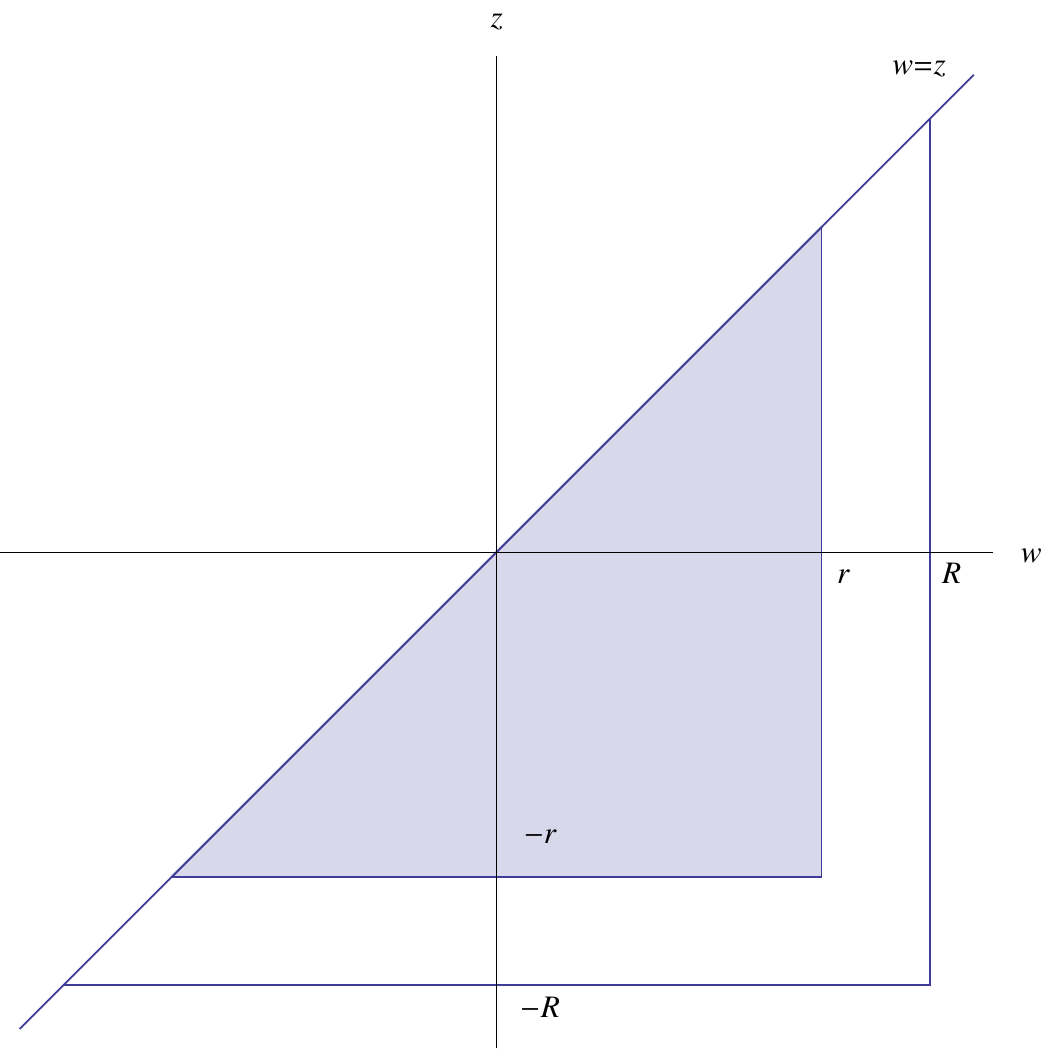}
  \caption{Assuming that the initial data are in the shaded region, we show 
  the existence of an $R$ such that the solution remains in the larger triangle. 
  The vacuum line is $w=z$.} 
  \label{fig:trekant1}
\end{figure}
We assume for the moment that $w^h,z^h$ satisfies an a priori bound of the form
\begin{equation}\label{e11}
w^h(x,t)\le R,\qquad z^h(x,t)\ge -R,\qquad (x,t)\in\R\X[0,\infty),
\end{equation}
for some constants $R>r$, and we will find a condition relating $r$ and $R$ 
under which \eqref{e11} can be justified. 

We first observe that if \eqref{e11} holds, then, for any 
$(x_1,t_1),(x_2,t_2)\in \R\X[0,\infty)$,
\begin{equation}\label{eL1}
|y^h(x_1,t_1)-y^h(x_2,t_2)|\le C(R)(|x_1-x_2|+|t_1-t_2|+h),  
\end{equation}
for some constant $C(R)>0$ depending only on $R$. 
In what follows, $C(R)$ will always represent a positive 
constant depending on $R$ that may differ from one occurrence to the next one. 

We also observe that 
\begin{equation}\label{eL2}
\begin{aligned}
&|\s^h(x_1,t_1)-\s^h(x_2,t_2)|\\
&\qquad=
|\tilde\s(y^h(x_1,t_1),t_1)-\tilde\s(y^h(x_2,t_2),t_2)|\\
&\qquad\le |\tilde\s(y^h(x_1,t_1),t_1)-\tilde\s(y^h(x_2,t_2),t_1)|
\\ & \qquad\qquad 
+|\tilde\s(y^h(x_2,t_2),t_1)-\tilde\s(y^h(x_2,t_2),t_2)|\\
&\qquad\le |y^h(x_1,t_1)-y^h(x_2,t_2)|\int_0^1|\tilde\s_y(\tau y_2^h
+(1-\tau)y_1^h,t_1)|\,d\tau\\
&\qquad\qquad+|t_1-t_2|\int_0^1|\tilde\s_t(y_2^h,\theta t_2
+(1-\theta)t_1)|\,d\theta\\
&\qquad\le C(R)\left( (|x_1-x_2|+|t_1-t_2|+h)e^{-t_1}
+|t_1-t_2|e^{-\min(t_1,t_2)}\right)\\
&\qquad\le C(R)(|x_1-x_2|+|t_1-t_2|+h) e^{-\min(t_1,t_2)},
\end{aligned}
\end{equation}
where we have used \eqref{e6''} and denoted $y^h_i=y^h(x_i,t_i)$, $i=1,2$.

\begin{figure}[tbp]
  \centering
  \includegraphics[width=0.65\linewidth]{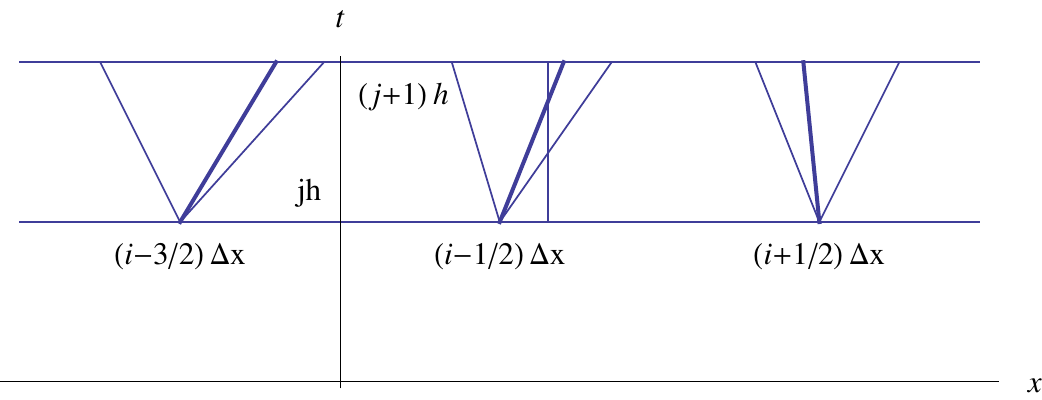} 
   \caption{Schematic figure of the solution of the Riemann problem. 
   Contact discontinuities are indicated by thick lines. 
   We see that the vertical line $x$ equals a constant 
   first intersects a slow or a fast wave before it crosses the contact discontinuity.} 
  \label{fig:waveB}
\end{figure}

Assume inductively that 
\begin{equation*}
w^h(x,t)\le r_j,\qquad    z^h(x,t) \ge-r_j, \qquad (x,t)\in\R\X[0,jh],
\end{equation*}
for some constant $r_j$. 
For $t\in[jh,(j+1)h)$ the approximate solution is defined 
by solving the Riemann problems given by the discontinuities 
at the points  $(x_{i+1/2},jh)$, $i\in\Z$. Since the $p$-system 
enjoys an invariant region given in terms 
of $w$ and $z$, the only possible increase in $w$ beyond $r_j$, 
and, similarly, the only possible decrease in $z$ beyond  $-r_j$, may 
occur across the contact discontinuity.  
Here both the velocity and the pressure remain 
unchanged, and the sole change is in the entropy.  
Observe first that since the slow Riemann invariant is 
increasing in $w$, there can be no increase in the value of $w$. 
Fix $x$ and let $t\in[jh,(j+1)h)$. We see from Figure \ref{fig:waveB} that 
the vertical line $x$ equals a constant crosses slow or 
fast waves before it crosses the contact discontinuity. 
Let $jh< \tilde t < \bar t <(j+1)h$ denote two times such that $\tilde t$ is 
after the fast or slow wave, but prior to the contact 
discontinuity, while $\bar t$  is after the contact discontinuity. Then we find  
\begin{equation*}
\begin{aligned}
z^h(x,\bar t)&= z^h(x,\tilde t)+(z^h(x,\bar t)-z^h(x,\tilde t)) \\ & 
\ge z^h(x,\tilde t) -\abs{z^h(x,\bar t)-z^h(x,\tilde t)}\\
&\ge -r_j-\abs{z^h(x,\bar t)-z^h(x,\tilde t)},
\end{aligned}
\end{equation*}
since the solution of  the $p$-system remains within the invariant region.
Furthermore,
\begin{equation}
\begin{aligned}
&\abs{z^h(x,\bar t)-z^h(x,\tilde t)}\\
&\quad\le 
\abs{(\frac{p^h}{\k}e^{-S^h/\Rf})^{\te/\g}(x,\bar t)
-(\frac{p^h}{\k}e^{-S^h/\Rf})^{\te/\g}(x,\tilde t)} \\
&\quad\le (\frac{p^h}{\k}e^{-S^h/\Rf})^{\te/\g}(x,\bar t)\, e^{\te S^h/(\g \Rf)}(x,\bar t) \\
&\qquad\qquad\times 
\abs{e^{-\te S^h/(\g \Rf)}(x,\tilde t)-e^{-\te S^h/(\g \Rf)}(x,\bar t)} \\
&\quad\le\frac12 (w^h-z^h)(x,\bar t)\, e^{\te S^h/(\g \Rf)}(x,\bar t)\\
&\qquad\qquad \times 
\abs{e^{-\te S^h/(\g \Rf)}(x,\tilde t)-e^{-\te S^h/(\g \Rf)}(x,\bar t)} \\
&\quad\le   r_j C(R) (\g-1)\abs{\jmp{S^h(\bar t)}}\\
&\quad =  r_j C(R) (\g-1)\abs{\jmp{S^h(jh)}},
\end{aligned}
\end{equation}
as we have replaced both the jump in the exponential by the 
corresponding jump in the exponent and estimated 
$e^{\te S^h/(\g \Rf)}(x,\bar t)$ by a common constant $C(R)$.  
Next we estimate the jump in the entropy. Let $x_1$ and $x_2$ be 
two points on the left and right side of a jump, respectively, thus 
$x_1<x_{i-1/2}<x_2$, with $x_2-x_1<\Dx$. We obtain  
\begin{equation} \label{eq:entroEst}
\begin{aligned}
\abs{\jmp{S^h(jh)}}&=\abs{S^h(x_2,jh)-S^h(x_1,jh)} \\
&\le \frac1{\Dx}\int_{I_i}|\s^h(\tilde x+\Dx,jh)
-\s^h(\tilde x,jh)|\, d\tilde x \\
&\le C(R) h\, e^{-jh},
\end{aligned}
\end{equation}
by \eqref{eL2}. 

This yields
\begin{equation}
z^h(x,\bar t)\ge -r_j(1+C(R) (\g-1)h\, e^{-jh}), 
\end{equation}
and we conclude that 
\begin{equation}
z^h(x,t)\ge -r_j(1+C(R) (\g-1)h\, e^{-jh}), \quad t\in [jh,(j+1)h).\label{eB1}
\end{equation}
A similar calculation leads to
\begin{equation}
w^h(x,t)\le r_j(1+C(R) (\g-1)h\, e^{-jh}), \quad t\in [jh,(j+1)h).\label{eB2}
\end{equation}

At $t=(j+1)h$ we average the approximate solution as described 
in \eqref{eq:snittB}--\eqref{eq:snitt}. Here we argue as follows. 
We first observe that the averaging of the values of 
$(\rho^h(x,(j+1)h-0), m^h(x,(j+1)h-0))$ in the intervals 
$I_i^{j+1}:=I_i\X\{t=(j+1)h-0\}$, $i\in\Z$, in order to obtain the 
values of $(\rho^h(x,(j+1)h), m^h(x,(j+1)h))$ in these intervals, does 
not affect the bounds \eqref{eB1} and \eqref{eB2}. More 
precisely, at each such interval, $S^h(x,(j+1)h-0)$ 
assumes at most 3 values, due to the possibility that two contact 
discontinuities, departing from $(x_{i-1/2},jh)$ and $(x_{i+1/2},jh)$, 
respectively, end inside $I_i^{j+1}$. This means that the 
values of $(\rho^h,m^h)$ in each interval $I_i^{j+1}$ belong to the 
union of at most 3 regions of the form
$$
R_\a:=\{(\rho,m)\,:\,-C\rho+e^{\theta S_\a/\Rf}\rho^{\theta+1}\le m
\le C\rho- e^{\theta S_\a/\Rf}\rho^{\theta+1}\},\qquad \a=1,2,3,
$$
for some constant $C>0$ common to all regions $R_\a$, $\a=1,2,3$. 
But, one easily check that $S_1<S_2$ implies $R_1\supset R_2$, 
that is, the regions $R_\a$, $\a=1,2,3$, are 
contained in that one corresponding to $S_*=\min\{S_1,S_2,S_3\}$. In 
particular, if we define
$$
S_*^h(x,(j+1)h):=\min\{S^h(\xi,(j+1)h-0)\,:\, \xi\in I_i\}, \qquad \text{in $I_i^{j+1}$, $i\in\Z$},
$$
then, from the convexity of the regions $R_\a$, we have 
\begin{equation}
z(\rho^h,u^h,S_*^h)(x,(j+1)h)\ge -r_j(1+C(R) (\g-1)h\, e^{-jh}), \label{eB1'}
\end{equation}
and also
\begin{equation}
w(\rho^h,u^h,S_*^h)(x,(j+1)h)\le r_j(1+C(R) (\g-1)h\, e^{-jh}), \label{eB2'}
\end{equation}
where 
$$
u^h(x,(j+1)h):=\begin{cases}\frac{m^h(x,(j+1)h)}{\rho^h(x,(j+1)h)},& 
\text{if $\rho^h(x,(j+1)h)>0$}\\ u^h(x,(j+1)h-0),&\text{otherwise}\end{cases} 
$$
and we agree that the value of $u^h(x,(j+1)h-0)$ at a 
vacuum interval is the mean value between its values at the 
extremes of the interval, which determines precisely 
the values of $u^h(x,(j+1)h-0)$ for all $x\in\R$. Observe 
also that the case in which $I_i^{j+1}$ is contained in a 
vacuum interval is trivial since $\rho^h=m^h=0$ in such an interval, and so the 
values of $\rho^h$ and $m^h$  do not change through averaging on $I_i^{j+1}$.

Now, we need to check how the bounds \eqref{eB1'} and \eqref{eB2'} 
change when we replace $S_*^h(x,(j+1)h)$ by the values of 
$S^h(x,(j+1)h)$ given by \eqref{eq:snitt}.
For this, we first estimate the change in $S^h$ from $S^h(x,(j+1)h-0)$, to
$S^h(x,(j+1)h)$, given by  \eqref{eq:snitt}. As already 
mentioned, $S^h(x,(j+1)h-0)$ can be one of three 
values; either the value $S^h(x,jh)$, or the values of $S$ in 
the neighboring intervals, that is, $S^h(x\pm\Dx,jh)$. In any of 
the three cases, the entropy is given 
by a formula similar to \eqref{eq:snitt}, but with $(j+1)h$ 
replaced by $jh$.  We consider the most representative case 
where the value is in a neighboring interval.  Thus
\begin{equation}
\begin{aligned}
& \abs{S^h(x,(j+1)h)-S^h(x-\Dx,jh)}
\\ & \qquad \le\frac{1}{\Delta x}\int_{I_i}
\abs{\s^h(\tilde x,(j+1)h)-\s^h(\tilde x-\Dx,jh)}\,d \tilde x \\
&\qquad \le C(R) h\, e^{-jh},
\end{aligned}
\end{equation}
again by \eqref{eL2}. Since,
$$
z^h(x,(j+1)h)=z(\rho^h,u^h,S^h)(x,(j+1)h),\ w^h(x,(j+1)h)=z(\rho^h,u^h,S^h)(x,(j+1)h),
$$
we conclude as above that 
\begin{equation}\label{edefrj}
\begin{aligned}
z^h(x,t)&\ge z(\rho^h,u^h,S_*^h)(x,(j+1)h)\\&
-|z(\rho^h,u^h,S_*^h)(x,(j+1)h)-z(\rho^h,u^h,S^h)(x,(j+1)h)|\\
&\ge -r_j(1+C(R) (\g-1)h\, e^{-jh})^2=: -r_{j+1}, \\
w^h(x,t)&\le w(\rho^h,u^h,S_*^h)(x,(j+1)h)\\&
+|w(\rho^h,u^h,S_*^h)(x,(j+1)h)-w(\rho^h,u^h,S^h)(x,(j+1)h)|\\
&\le r_j(1+C(R) (\g-1)h\, e^{-jh})^2 = r_{j+1}.
\end{aligned}
\end{equation} 

It remains to estimate the $r_j$. {}From the inductive 
formula \eqref{edefrj} for the $r_j$, we find
\begin{equation}\label{edefrj2}
\begin{aligned}
r_j&=r\prod_{k=1}^j(1+C(R) (\g-1)h\, e^{-kh})^2  \\
&\le r\exp\big( 2C(R)(\g-1)\sum_{k=1}^j e^{-kh}h)  \big)\\
&\le r\exp\big(2C(R) (\g-1)\int_{0}^\infty e^{-s}ds  \big)\\
&\le r e^{2C(R) (\g-1)}.
\end{aligned}
\end{equation}

Therefore, we see from \eqref{edefrj2} that the condition relating $r$ and $R$ under which the a priori bound \eqref{e11} holds is 
\begin{equation}
R e^{-2(\g-1)C(R)}\ge r.\label{e12}
\end{equation}
We may easily check that $C(R)$ may be defined as a continuous increasing function of $R\in[0,\infty)$ such that $C(0)=0$ and $C(R)\to\infty$ as $R\to\infty$.
Hence, the left-hand side of \eqref{e12} attains a maximum value for some $R_*\in(0,\infty)$ and by \eqref{e12} the initial bound $r$ can take the largest possible value given by the left-hand side of \eqref{e12} for $R=R_*$. In particular, \eqref{e12} may be viewed as a restriction on the initial bound $r$ which amounts a restriction on $\|\rho_0\|_\infty$ and $\|m_0\|_\infty$, assuming given $S_0$. We also verify that the initial bound can be taken as large as we wish provided that  $\g-1$ is sufficiently small.

\medskip
Now we proceed to prove the compactness of the sequence of approximate solutions $U^h$. The proof is based on the general analysis carried out by DiPerna in 
\cite{DP} and we are going to apply the compactness result in \cite{DP2} and its extensions in \cite{Ch}, \cite{LPT} and \cite{LPS}, which together cover the whole range $\g>1$.

Now, let $V^h=(\rho^h,m^h)$ and $F^h=(m^h,\rho^h(u^h)^2+p(\rho^h,S^h))$. For any $\phi\in C_0^\infty(\R^2)$ we have  
\begin{equation}\label{e13}
\begin{aligned}
\iint_{\R\X [0,\infty)}V^h\phi_t+F^h\phi_x\,dx\,dt &=
\sum_{j=0}^\infty\int_{jh}^{(j+1)h}\int_{\R}V^h\phi_t+F^h\phi_x\,dx\, dt \\
&=\sum_{j=0}^\infty\int_\R\jmp{V^h(jh)}\phi(x,jh)\,dx \\
&=\sum_{j=1}^\infty\int_\R\jmp{V^h(jh)}\phi(x,jh)\,dx
-\int_\R V^h(x,0)\phi(x,0)\,dx,
\end{aligned}
\end{equation}
where 
$$
\jmp{V^h(jh)}=V^h(x,jh-0)-V^h(x,jh+0). 
$$
Further, if $(\eta,q)$ is an arbitrary entropy pair for \eqref{e1}--\eqref{e2}, with $S$ constant, we have
\begin{equation}\label{e14}
\begin{aligned}
&\iint_{\R\X [0,\infty)}\eta^h\phi_t+q^h\phi_x\,dx\,dt=
\sum_{j=0}^\infty \int_{jh}^{(j+1)h}\int_{\R}\eta^h\phi_t+q^h\phi_x\,dx\,dt\\
&=-\int_\R\eta^h(x,0)\phi(x,0)\,dx
+\sum_{j=1}^\infty\int_\R\jmp{\eta^h(jh)}\phi(x,jh)\,dx 
+\int_0^\infty\CS(\phi)\,dt+\int_0^\infty\CC(\phi)\,dt,
\end{aligned}
\end{equation}
where, for reasons of brevity, we write $\eta^h=\eta(V^h,S^h)$ and $q^h=q(V^h,S^h)$. 
Here
\begin{equation*}
\jmp{\eta^h(jh)}=\eta^h(x,jh-0)-\eta^h(x,jh+0),
\end{equation*}
and $\CS(\phi)$ is defined as
\begin{align*}
&\CS(\phi)=\sum_{\text{shocks}}\big(s\jmp{\eta^h}-\jmp{q^h}\big)\phi(x(t),t),\\
&\jmp{\eta^h}=\eta^h(x(t)-0,t)-\eta^h(x(t)+0,t),
\end{align*}
where the sum is over all shock discontinuities $(x(t),t)$ at time $t$, $s=x'(t)$ denoting the shock speed, while $\CC(\phi)$ is defined as
$$
\CC(\phi)=\sum_{\substack{\text{contact}\\ \text{discontinuities}}}
\big(u^h\jmp{\eta^h}-\jmp{q^h}\big)\phi(x(t),t),
$$
with sum running over all contact discontinuities $(x(t),t)$ at time $t$,  where 
$u^h$ is the velocity. The latter is defined over a vacuum interval as the arithmetic mean between the velocity at the end of the 1-rarefaction wave bounding the 
vacuum interval on the left-hand side   and the velocity at the beginning of the 2-rarefaction wave 
bounding the vacuum interval on the right-hand side.   

We recall that if $(\eta,q)$ is a convex entropy pair for the isentropic 
system \eqref{e1}--\eqref{e2} where $S$ is constant, then 
\begin{equation}\label{eent2}
s\jmp{\eta^h}-\jmp{q^h}\ge 0,
\end{equation}
across each shock wave.  Since $S^h$ is constant across waves 
of the first and third family, inequality \eqref{eent2} also holds here.
Therefore, for any weak entropy pair $(\eta,q)$, we find that the functional
$$
\int_0^\infty \CS(\phi)\,dt
$$
is a (signed) measure with locally finite total variation, as a consequence of Remark~\ref{R:1}. 

Concerning the functional 
$$
\int_0^\infty \CC(\phi)\,dt,
$$
if $(\eta,q)$ is a smooth entropy pair, we have, in view of 
previous calculations, 
$$
|u^h\jmp{\eta^h(jh)}-\jmp{q^h(jh)}|\le C_\eta e^{-jh}h,
$$
and so
$$
\left|\int_0^\infty \CC(\phi)\,dt\right|\le C_\eta \operatorname{diam}(K)\|\phi\|_\infty,
$$
where $K$ is any compact containing the support of $\phi$, which  gives that this functional is also a measure with locally finite total variation. 

Observe that the weak entropies may be also written as
$$
\eta(\rho,u)=\rho\int_{-1}^1g\left(\frac{m}{\rho}+ze^{(\g-1)S/2R}\rho^{(\g-1)/2}\right)(1-z^2)_+^\l\,dz,
$$
while a similar formula holds for $q$. In particular, $\eta,q$ are Lipschitz up to vacuum if $g$ is smooth.  

We also observe that for the special entropy pair 
$(\eta_*,q_*)$ we have $\int_0^\infty \CC(\phi)\,dt=0$. 
Also, for this entropy pair, for nonnegative $\phi\in C_0^\infty(\R^2)$ we have 
\begin{equation}\label{eineq0}
\begin{aligned}
&\sum_{j=1}^\infty\int_\R\jmp{\eta_*^h(jh)}\phi(x,jh)\,dx\\
&\quad=\sum_{j=1}^\infty \sum_{i\in\Z}\int_{I_i} \Big(\eta_*(V^h(x,jh-0)),S^h(x,jh+0))\\
&\qquad\qquad\qquad\qquad\qquad\qquad
             -\eta_*(V^h(x,jh+0),S^h(x,jh+0))\Big)\phi(x,jh)\,dx\\
&\quad- \sum_{j=1}^\infty \sum_{i\in\Z}\int_{I_i}\Big(\eta_*(V^h(x,jh-0)),S^h(x,jh+0))\\
&\qquad\qquad\qquad\qquad\qquad\qquad-\eta_*(V^h(x,jh-0),S^h(x,jh-0))\Big)\phi(x,jh)\,dx.
\end{aligned}
\end{equation}
The first sum in the right-hand side of equation \eqref{eineq0} is nonnegative for nonnegative $\phi$, since $V^h(x,jh+0)$ is the average of $V^h(x,jh-0)$, in each interval $I_i$, and $\eta_*$ is convex. 
Therefore, we get
\begin{equation}\label{eineq1}
\begin{aligned}
& \sum_{j=1}^\infty\int_\R\jmp{\eta_*^h(jh)}\phi(x,jh)\,dx \\
&\quad\ge -\sum_{j=1}^\infty \sum_{i\in\Z}
\int_{I_i}\eta_{*S}^h(\cdots)(S^h(x,jh+0)-S^h(x,jh-0))\phi(x,jh)\,dx\\
&\quad\ge-\sum_{j=1}^\infty Ce^{-jh}h\int_{\R} \phi(x,jh)\,dx,
\end{aligned}
\end{equation}
where $\eta_{*S}^h(\cdots)=\int_0^1\eta_{*S}^h(V^h(x,jh-0),
A(\theta))\,d\theta$ is the coefficient of the linear remaining 
term in the trivial Taylor expansion of zero order in the variable $S$ and
$A(\theta)=(1-\theta)S^h(x,jh-0)+\theta S^h(x,jh+0)$. 
In particular, both the left-hand side as well as the second term of 
the right-hand side of \eqref{eineq0} are measures of locally finite total variation. 
As a consequence, we may apply equality \eqref{eineq0} with $\phi$ replaced by the 
characteristic function of any suitably chosen 
rectangle $|x|\le L=M\Delta x$, $0\le t\le T=Nh$, to find that
\begin{equation}
\sum_{jh\le N} \sum_{|i\Delta x|\le M}\int_{I_i}D_V^2
\eta_*^h(\cdots)(V^h(x,jh-0))-V^h(x,jh+0))^2\,dx\le \text{const.},\label{eineq2}
\end{equation}
for any $M,N>0$, the constant depending on $M,N$, where 
$D_V^2\eta_*^h(\cdots)=\int_0^1(1-\theta) 
D_V^2\eta_*(B(\theta),S^h(x,jh+0))\,d\theta$ is the coefficient 
of the  quadratic remaining term in the Taylor expansion of first order and
$B(\theta)=(1-\theta)V^h(x,jh+0)+\theta V^h(x,jh-0)$.  

Since for all weak entropy $\eta$ we have $|D_V^2\eta|\le C_\eta D_V^2\eta_*$, 
for some $C_\eta>0$, it follows from \eqref{eineq2} that
\begin{equation}\label{eineq3}
\left|\sum_{jh\le N} \sum_{|i\Delta x|\le M}
\int_{I_i}|D_V^2\eta|(V^h(x,jh-0))-V^h(x,jh+0))^2\,dx\right|\le \text{const.},
\end{equation}
for any $M,N>0$, the constant depending on $M,N$. 

We can then use DiPerna's method in \cite{DP} to 
prove the $W_\loc^{-1,2}$ compactness of
the distributions $\eta_t^h+q_x^h$ by decomposing the functional
$$
L(\phi)=\sum_{j=1}^\infty\int_\R\jmp{\eta^h(jh)}\phi(x,jh)\,dx
$$
as 
\begin{equation}
\begin{aligned}
L(\phi) 
&= \sum_{j=1}^\infty\sum_{i\in\Z}\int_{I_i}\jmp{\eta^h(jh)}\phi(x,jh)\,dx   \\
&=\sum_{j=1}^\infty\sum_{i\in\Z}\Big(\phi(x_i,jh)\int_{I_i}\jmp{\eta^h(jh)}\,dx\\
&\qquad\qquad\qquad
+ \int_{I_i}\jmp{\eta^h(jh)}\big(\phi(x,jh)-\phi(x_i,jh)\big)dx\Big)   \\
&=L_1(\phi)+L_2(\phi).
\end{aligned}
\end{equation}
We consider the two terms separately. We have
$$
L_1(\phi)=\sum_{j=1}^\infty\sum_{i\in\Z}\phi(x_i,jh)
\int_{I_i}\jmp{\eta^h(jh)}_V+\jmp{\eta^h(jh)}_S\,dx=:L_{11}(\phi)+L_{12}(\phi),
$$
where, if $\jmp{\eta(V,S)}=\eta(V_-,S_-)-\eta(V_+,S_+)$, we denote 
$$
\jmp{\eta(V,S)}_V=\eta(V_-,S_-)-\eta(V_+,S_-), \qquad 
\jmp{\eta(V,S)}_S=\eta(V_+,S_-)-\eta(V_+,S_+).
$$
Since $|\jmp{\eta^h(jh)}_S|\le Ce^{-jh}h$, we clearly have 
$$
|L_{12}(\phi)|\le C\|\phi\|_\infty.
$$
Concerning $L_{11}(\phi)$, we have, cf.~\eqref{eineq3}, 
\begin{equation}
\begin{aligned}
\abs{L_{11}(\phi)}&\le\left|\sum_{j=1}^\infty\sum_{i\in\Z}\phi(i,jh) 
\int_{I_i}\jmp{\eta^h(jh)}_V\,dx\right| \\
&=\left|\sum_{j=1}^\infty\sum_{i\in\Z}\phi(i,jh) 
\int_{I_i}D_V^2\eta^h(\cdots)(V^h(x,jh-0)-V^h(x,jh+0))^2\,dx\right|\\
&\le  C\norm{\phi}_\infty. 
\end{aligned}
\end{equation}
Hence, we have
$$
|L_1(\phi)|\le C_1\|\phi\|_\infty.
$$

Next, exactly as in \cite{DP}, we find, assuming 
that the test function $\phi$ satisfies 
$\supp{\phi}\subseteq [-N,N]\X[-J,J]$ and 
keeping $\theta>0$ sufficiently small, 
\begin{equation}
\begin{aligned}
\abs{L_2(\phi)}&\le\sum_{\abs{j}\le J}\sum_{\abs{i}\le N}
\int_{I_i}\abs{\jmp{\eta^h(jh)}}\abs{\phi(x,jh)-\phi(i,jh)}\,dx  \\
&\le \norm{\phi}_{C^\a} \sum_{\substack{\abs{j}\le J\\ \abs{i}\le N}}
\int_{I_i}\abs{\jmp{\eta^h(jh)}}\Dx^\a\,dx\\
&\le  \norm{\phi}_{C^\a} \sum_{\substack{\abs{j}\le J\\ \abs{i}\le N}}\int_{I_i}
\Big(\frac{\Dx^{2\a}}{\Dx^\te}+\Dx^\te\abs{\jmp{\eta^h(jh)}}^2 \Big)dx \\
&\le  \norm{\phi}_{C^\a} \sum_{\abs{j}\le J}\sum_{\abs{i}\le N}
\Big(\frac{\Dx^{2\a+1}}{\Dx^\te}+\Dx^\te\int_{I_i}\abs{\jmp{\eta^h(jh)}}^2 \, dx\Big) \\
&\le  \norm{\phi}_{C^\a}\Big(\frac{\Dx^{2\a+1}}{\Dx^\te}(2J+1)(2N+1)+\Dx^\te
\sum_{\abs{j}\le J}\int_{\R}\abs{\jmp{\eta^h(jh)}}^2 \, dx\Big) \\ 
&\le  \norm{\phi}_{C^\a}\Big(\frac{\Dx^{2\a+1}}{\Dx^\te}\Oh\big(\frac{1}{\Dx\Dt}\big)
+\Dx^\te\sum_{\abs{j}\le J}\int_{\R}\abs{\jmp{\eta^h(jh)}}^2 \, dx\Big)\\
&\le  C_2\norm{\phi}_{C^\a}\Big(\frac{\Dx^{2\a+1}}{\Dx^{\te+2}}+\Dx^\te\Big) \\ 
&\le C_2\norm{\phi}_{C^\a} \Dx^{\a-1/2}
\end{aligned}
\end{equation}
where $C^\a$ denotes the H\"older space with seminorm 
$$
\norm{\phi}_{C^\a}=\sup_{x,y\in\R}\abs{\phi(x)
-\phi(y)}/\abs{x-y}^\a, \qquad \a>1/2,
$$ 
and where $C_2$ depends on the support of 
$\phi$. Thus
$$
\abs{L_1(\phi)}\le C_1\|\phi\|_\infty,\qquad\text{and}\qquad 
\abs{L_2(\phi)}\le C_2(\Delta x)^\b\|\phi\|_{C^\a}
$$
for appropriate $\a,\b\in(0,1)$, for some positive constants $C_1,C_2$ 
depending on  $\supp \phi$, but independent of $\phi$, and 
through the Sobolev imbedding theorem
$$
L_2(\phi)\le C_2(\Delta x)^\b\|\phi\|_{W^{1,q}},
$$
for an appropriate $q\in (1,2)$ and constant depending on the support of $\phi$.  

In this way we obtain by the usual interpolation argument that for 
any weak entropy pair $(\eta,q)$ for \eqref{e1}--\eqref{e2} we have
$$
\eta(V^h,S^h)_t+q(V^h,S^h)_x\in\{\ \text{compact of } W_\loc^{-1,2}(\R\X[0,\infty))\ \}.
$$
We can then use the compactness results in \cite{DP2,Ch,LPT,LPS} to 
deduce that we may extract a subsequence of $(\rho^h,m^h,S^h)$ 
converging in $L_\loc^1(\R\X[0,\infty))$ to a weak solution 
$(\rho(x,t),m(x,t),S(x,t))$ to \eqref{e1}--\eqref{e5}. Also, \eqref{eineq1} 
implies the entropy inequality \eqref{eent}, and \eqref{e13} 
implies \eqref{e8} by a calculation 
similar to the estimate for $L_2(\phi)$ above.  

\bigskip
Concerning the decay property \eqref{edecay}, we prove it as follows. 
First, from the above discussion, we deduce that for any weak entropy pair we have
$$
|\la \eta(\rho,m,S)_t+q(\rho,m,S)_x,\phi\ra|\le C_1\|\phi\|_\infty,
$$
with $C_1$ depending only on $\supp \phi$ and bounds for $(\rho,m,S)$. 
Hence, if $U^T=(\rho^T,m^T,S^T)$ is the self-scaling 
sequence $U^T(x,t)=U(Tx,Tt)$, we see that for any entropy  pair
$$
|\la \eta(\rho^T,m^T,S^T)_t+q(\rho^T,m^T,S^T)_x,\phi\ra|\le C_1\|\phi\|_\infty,
$$
while from \eqref{eent} we have, for $0\le t\le T$, 
\begin{align*}
 \int_{[0,L]} \eta_*(\rho,m,S)(x,t)\,dx &\ge \int_{[0,L]} \eta_*(\rho,m,S)(x,T)\,dx 
 -C\int_t^T\int_{[0,L]}e^{-s}\,dx\,ds\\
  &\ge\int_{[0,L]} \eta_*(\rho,m,S)(x,T)\,dx- CLe^{-t}
\end{align*}
Hence, we can apply the decay result in \cite{CF} to 
deduce \eqref{edecay}, which then concludes the proof. 


\end{document}